\documentclass[12pt, reqno, letterpaper]{amsart}
\usepackage{amssymb}  
\usepackage{bbold}  

\theoremstyle{plain}
\newtheorem{theorem}{Theorem} 
\newtheorem{definition}[theorem]{Definition}
\newtheorem{lemma}[theorem]{Lemma}
\theoremstyle{remark}
\newtheorem*{acknowledgement}{Acknowledgement}

\newcommand\res{\mathord {\upharpoonright}}  

\newcommand\BB{{\mathcal B}}
\newcommand\CC{{\mathcal C}}
\newcommand\PP{{\mathcal P}}

\newcommand{\PPP}{\mathbb{P}}

\newcommand\hgt{\mathrm{ht}}
\newcommand\rank{\mathrm{rank}}  
\newcommand\Mrank{\mathrm{Mrank}}  

\newcommand\one{\mathbb{1}} 

\newcommand\supt{\mathrm{supt}}  

\title{Compact Scattered Spaces in Forcing Extensions}

\author{Kenneth Kunen}

\address{Department of Mathematics \\
University of Wisconsin \\
Madison, WI 57306 USA}

\email{kunen@math.wisc.edu}

\urladdr{http://www.math.wisc.edu/\symbol{126}kunen/}

\date{\today}



\thanks{Author partially supported by NSF Grant DMS-0097881.}

\begin{document}

\begin{abstract}
We consider the cardinal sequences of compact
scattered spaces in models where CH is false.
We describe a number of models of $2^{\aleph_0} = \aleph_2$
in which no such space can have $\aleph_2$ countable levels.
\end{abstract}

\maketitle

All spaces considered here are assumed to be Hausdorff.
The following definition summarizes the standard 
Cantor-Bendixson sequence:

\begin{definition}
For any topological space $X$:
\begin{itemize}
\item[1.] $I(X)$ is the set of isolated points of $X$.
\item[2.] $X^{(0)} = X$. 
\item[3.] $X^{(\alpha + 1)} = X^{(\alpha)} \setminus I(X^{(\alpha)})$.
\item[4.] $X^{(\gamma)} = \bigcap_{\alpha < \gamma} X^{(\alpha)}$ for limit
$\gamma$.
\item[5.] $X$ is \emph{scattered}
iff $X^{(\alpha)} = \emptyset$ for some $\alpha$; then, the least such 
$\alpha$ is called $\hgt(X)$, the \emph{height} of $X$.
\end{itemize}
\end{definition}
\noindent
Equivalently, $X$ is scattered iff $I(Y) \ne \emptyset$ for all
nonempty $Y \subseteq X$.
If $X$ is compact scattered and nonempty, then $\hgt(X)$ is a successor
ordinal, $\delta+1$, and $X^{(\delta)}$ is finite.

Juh\'asz, Shelah, Soukup, and Szentmikl\'ossy \cite{JSSS} study the possible
values for the \textit{cardinal sequence},
$\langle |I(X^{(\alpha)})| : \alpha < \hgt(X) \rangle $,
for scattered spaces $X$.  \cite{JSSS}, combined with earlier work,
shows that the class of cardinal sequences obtained from the
regular scattered spaces is determined by cardinal arithmetic,
but this is not true for compact spaces.
In particular,  consider the following assertions:
\begin{itemize}
\item[$(A)$] There is a compact scattered $X$
such that $\hgt(X) = \omega_2+1$ and $|I(X^{(\alpha)})| = \aleph_0$ for
all $\alpha < \omega_2$. 
\item[$(B)$]
$|\{\alpha < \hgt(X) : |I(X^{(\alpha)})| = \aleph_0\}| \le \aleph_1$ for
every compact scattered $X$.
\end{itemize}
Clearly $(B)$ implies $\neg (A)$.
Baumgartner and Shelah \cite{BS} showed that $(A)$ is consistent
with $2^{\aleph_0} = \aleph_2$,
whereas Just \cite{Ju} showed that $\neg (A)$ is consistent
with $2^{\aleph_0} = \aleph_2$.
Just's model was the standard Cohen real extension over a model of CH.
This result was improved by \cite{JSSS},  which showed that
$(B)$ holds in the Cohen extension.
The argument in \cite{JSSS} is very specific to Cohen forcing,
and they ask whether $(B)$ holds in similar extensions,
such as the one obtained by adding $\aleph_2$ random reals side-by-side.
We show here (see Theorem \ref{thm-main})
that in fact $(B)$ does hold in this model,
and also in the model obtained by adding $\aleph_2$ random reals
in the standard way (by one measure algebra).
We do not know whether $(B)$ follows from some abstract principle,
such as $C^s(\omega_2)$ described by 
Juh\'asz,  Soukup, and  Szentmikl\'ossy \cite{JSS-comb};
this is also asked in \cite{JSSS}.
However, since $(B)$ is true in the random real model,
where $C^s(\omega_2)$ is false (because there is an $\omega_2$--Luzin gap),
it would be more interesting if one could derive $(B)$ from
some abstract principle true in both the Cohen and random extensions.

Some remarks: \cite{JSSS} states $(B)$ for locally compact spaces,
but that is equivalent (by taking the one-point compactification).
Likewise, if $(A)$ holds, one can remove the top level and
get a locally compact $X$ of height $\omega_2$ such that
all $I(X^{(\alpha)})$ are countable.
Under CH, $(B)$ is true, since if $\alpha$ is least with
$I(X^{(\alpha)})$ countable, then the weight of
$X^{(\alpha)}$ is no more than  $2^{\aleph_0} = \aleph_1$
so that $\hgt(X) < \alpha + \omega_2$.

Similarly to \cite{JSSS}, our argument assumes that $(B)$ is false,
and produces an infinite independent sequence of clopen sets,
contradicting the assumption that $X$ is scattered.
We begin with some elementary remarks on such sequences:

\begin{definition}
Assume that $\nu \le \omega$ and $K_i \subseteq X$ for $i < \nu$.
Let $K_i^0 = K_i$ and
$K_i^1 = X \backslash K_i$.  The $\nu$-sequence
$\langle K_i : i < \nu \rangle$ is \emph{independent} iff
$\bigcap_{i < n} K_i^{s(i)} \ne \emptyset$ for each \emph{finite}
$n \le \nu$ and each $s \in \{0,1\}^n$.
\end{definition}

\begin{lemma}
\label{lemma-ind-clopen}
If $X$ is compact scattered, then there is no independent $\omega$-sequence
of clopen subsets of $X$.
\end{lemma}

If $X$ is also infinite, there will be independent
$n$-sequences of clopen subsets for each finite $n$.
These finite sequences form a tree in the natural way,
ordered by extension; the root of the tree
is the empty sequence, $\emptyset$.
Then by Lemma \ref{lemma-ind-clopen},
this tree is well-founded --- that is,
it can have no infinite paths.
Our proof of Theorem \ref{thm-main}
will require a somewhat more complicated tree:

\begin{definition}
\label{def-tree}
Suppose that $\CC_\xi \subseteq \PP(D)$ for each $\xi < \gamma$.
Define the tree $T = T( \langle \CC_\xi : \xi < \gamma \rangle)$
as follows:
\begin{itemize}
\item[1.] The nodes of $T$ at level $n$ are pairs $(\vec \xi, \vec H)$, where
$\vec \xi = \langle \xi_i : i < n \rangle\allowbreak\in \gamma^n$
is a sequence of distinct ordinals and 
$\vec H = \langle H_i : i < n \rangle \in (\PP(D))^n$
is an independent sequence, with
each $H_i \in \CC_{\xi_i}$.
\item[2.] $<$ denotes the usual tree order, with the root
$(\emptyset,\emptyset)$ at the top.
\end{itemize}
If $T$ is well-founded, let $\rho = \rho_T$ be its
\emph{rank function}, defined by
$\rho(x) = \sup\{\rho(y) + 1 : y \in T \ \&\ y < x\}$;
also, define $\rank(T) = \rho_T( \emptyset,\emptyset )$,
and define $\Mrank(\langle \CC_\xi : \xi < \gamma \rangle)$
to be the minimum possible value among all
$\rank(T(\langle \CC_{\xi_\mu} : \mu < \gamma \rangle))$,
where $\langle \xi_\mu : \mu < \gamma \rangle$ is a strictly increasing
$\gamma$--sequence of ordinals less than $\gamma$.
\end{definition}

Thus, leaf nodes (ones with no proper extension in $T$) will
have rank $0$, and the root $(\emptyset, \emptyset)$ will have the largest rank,
which we are calling $\rank(T)$. 
Then $T(\langle \CC_{\xi_\mu} : \mu < \gamma \rangle)$ will 
always be a subtree of $T$, so its rank is $\le \rank(T)$;
and $\Mrank(\langle \CC_\xi : \xi < \gamma \rangle)$
is the least among the ranks of these subtrees.
Note that Definition \ref{def-tree} does not refer to any topology
on $D$.  The topology arises in the following lemma:

\begin{lemma}
\label{lemma-bigrank}
Let $\gamma$ be any limit ordinal.
Assume that:
\begin{itemize}
\item[1.]  $X$ is compact scattered.
\item[2.]   $E_\xi$, for $\xi < \gamma$, are disjoint nonempty
subsets of $X$, and
$\overline{E_\eta} \supseteq E_\xi$ 
whenever $\eta < \xi < \gamma$.
\item[3.] Each $\BB_\xi$ is a subalgebra of the clopen
subsets of $X$ which separates the points of $E_\xi$.
\end{itemize}
Then $\Mrank(\langle \BB_\xi : \xi < \gamma \rangle) \ge \gamma$
\end{lemma}
\begin{proof}
Let $T =  T( \langle \BB_\xi : \xi < \gamma \rangle)$.
It is sufficient to prove that $\rank(T) \ge \gamma$
because every subsequence of $\langle E_\xi : \xi < \gamma \rangle$
has the same properties.

Call $x = (\vec \xi, \vec H)$ at level $n > 0$
of $T$ \textit{special}
iff $\xi_0 > \xi_1 > \cdots > \xi_{n-1}$ and
$E_{\xi_{n-1}} \cap \bigcap_{i < n} H_i^{s(i)} \ne \emptyset$
for each  $s \in 2^n$.
We prove that 
$\rho(x) \ge \xi_{n-1}$ for such special nodes $x$.
That will imply the desired lemma, since such special nodes
clearly exist (even with $n=1$) for each possible value of
$\xi_{n-1} < \gamma$.

The proof proceeds by induction on $\xi_{n-1}$.
Thus, fix $x = (\vec \xi, \vec H)$.
It is sufficient to prove that for each $\eta < \xi_{n-1}$:
$x$ has an extension of the form $y = (\vec \xi', \vec H')$,
where $y$ is special, 
$\vec \xi' = \langle \xi_0, \ldots , \xi_{n-1}, \xi_n \rangle$,
$\xi_n = \eta$, and $\vec H'$ extends $\vec H$.
Then $H'$ will be of the form
$\langle H_0, \ldots , H_{n-1}, H_n \rangle$,
and, we need to define $H_n$.
Since each clopen set $\bigcap_{i < n} H_i^{s(i)}$  meets
$E_{\xi_{n-1}}$, and 
$\overline{E_\eta} \supseteq E_{\xi_{n-1}}$,
each $\bigcap_{i < n} H_i^{s(i)}$  meets $E_\eta$ in an infinite set.
Since $\BB_\eta$ is an algebra and
separates the points of $E_\eta$, we can choose $H_n$ so that it 
and its complement meet all the $2^n$
sets $E_\eta \cap \bigcap_{i < n} H_i^{s(i)}$.
\end{proof}

We remark that if $\hgt(X) \ge \gamma$,
then (1)(2) are satisfied by taking $E_\xi = I(X^{(\xi)})$.
If all the $E_\xi$ of Lemma \ref{lemma-bigrank} are countable,
then we obtain the situation of Definition \ref{def-tree}
with $D = \omega$:

\begin{lemma}
\label{lemma-omega}
Suppose that $X$, $\gamma$, and the $E_\xi$ satisfy
$(1)(2)$
from Lemma  \ref{lemma-bigrank}, and in addition all $|E_\xi| = \aleph_0$.
Then there are countable subalgebras
$\CC_\xi \subseteq \PP(\omega)$ for $\xi < \gamma$ such that
$T( \langle \CC_\xi : \xi < \gamma \rangle)$ is well-founded and
$\Mrank(\langle \CC_\xi : \xi < \gamma \rangle) \ge \gamma$.
\end{lemma}
\begin{proof}
WLOG, $X = \overline{E_0}$ (if not, replace $X$ by $\overline{E_0}$).
Then, WLOG, $E_0 = \omega$.
Choose $\BB_\xi $ as in (3) of Lemma \ref{lemma-bigrank},
with $\BB_\xi$ countable.
Let $\CC_\xi = \{H \cap \omega : H \in \BB_\xi\}$.
Then $\CC_\xi$ is a countable subalgebra of $\PP(\omega)$.
Observe that the map $H \mapsto H \cap \omega$ is an isomorphism
from $\BB_\xi$ onto $\CC_\xi$; its inverse is the map
$K \mapsto \overline K$.  Thus, 
$T( \langle \BB_\xi : \xi < \kappa \rangle)$
of Definition \ref{def-tree} is isomorphic to
$T = T( \langle \CC_\xi : \xi < \kappa \rangle)$,
so $\Mrank(T) \ge \gamma$ by Lemma \ref{lemma-bigrank}.
\end{proof}

We now turn to forcing extensions.  As usual (see, e.g., \cite{Jech, KUN}),
a \textit{partial order}
$\PPP$ really denotes a triple, $(\PPP, \le, \one)$, where
$\le$ is a transitive reflexive relation on $\PPP$ and $\one $ is
a largest element of $\PPP$.  Then,
$\prod_{i \in \theta} \PPP_i$ denotes the product of the
$\PPP_i$, with the natural product order.
Elements $\vec p \in \prod_{i \in \theta} \PPP_i$
are $\theta$-sequences, with each $p_i \in \PPP_i$.
The \textit{finite support product} is given by:

\begin{definition}
\label{def-fin-spt}
If $\vec p \in \prod_{i \in \theta} \PPP_i$,
then the {\em support} of $\vec p$, $\supt(\vec p)$, is
$\{i \in \theta : p_i\nobreak\ne\nobreak\one \}$.
$ \prod^{fin}_{i \in \theta} \PPP_i =
\{\vec p \in \prod_{i \in \theta} \PPP_i :
|\supt(\vec p)| < \aleph_0\}$.
\end{definition}

If all $\PPP_i$ are countable and non-atomic, one gets the
\textit{Cohen real extension}, adding $\theta$ Cohen reals
(note that this is the same extension for $1 \le \theta \le \aleph_0)$.
To get the \textit{random real extension}
(see \cite{Jech}), $\PPP$ is a measure algebra;
for example, the measure algebra of $\{0,1\}^\theta$ adds $\theta$ 
random reals when $\theta$ is infinite.  This is different from the 
\textit{side-by-side random real extension}, which uses
$ \prod^{fin}_{i \in \theta} \PPP_i$, where each $\PPP_i$ is
the measure algebra of $\{0,1\}^{\omega}$.
Our main result applies to both the random real extension
and the side-by-side random real extension.
First, a lemma about ranks whose conclusion goes in the 
opposite direction from that of Lemma \ref{lemma-omega}:

\begin{lemma}
\label{lemma-smallrank}
In the ground model $V$, set $\kappa = (2^{\aleph_0})^+$.
Let $\PPP$ be ccc and be  either a random real extension
\textup{(}adding any number of random reals\textup{)},
or of the form $\prod^{fin}_{i \in \theta} \PPP_i$,
where all the $\PPP_i$ are isomorphic and $|\PPP_i| \le 2^{\aleph_0}$
\textup{(}but $\theta$ is arbitrary\textup{)}.

Then in the generic extension $V[G]$, the following holds:
If
$\langle \CC_\xi : \xi < \kappa \rangle$ is 
as in Definition \ref{def-tree}, with $D = \omega$, all
$|\CC_\xi| = \aleph_0$, and
$T( \langle \CC_\xi : \xi < \kappa \rangle)$ well-founded, then 
$\Mrank(\langle \CC_\xi : \xi < \kappa \rangle) < \omega_1$.
\end{lemma}
\begin{proof}
We actually obtain the appropriate subsequence
$\langle \xi_\mu : \mu < \kappa \rangle$ in $V$
by a standard thinning-out process.
We first consider the case that
$\PPP = \prod^{fin}_{i \in \theta} \PPP_i$,
and then comment on what needs to be changed if $\PPP$ is a random
real forcing.

In $V$, for $J \subseteq \theta$,
let $\PPP_J = \prod^{fin}_{j \in J} \PPP_j$,
which we identify as a suborder of $\PPP$.
$\PPP_\emptyset$ denotes the one-element order $\{\one\}$.
We have, for each $\xi < \kappa$, a name
$\dot \CC_\xi$ which is forced by $\one$
to denote a countable subset of $\PP(\omega)$.
Thus, $\dot \CC_\xi$  is really a $\PPP_{J_\xi}$-name, where
$J_\xi \subseteq \theta$ is countable.
Since $\kappa = (2^{\aleph_0})^+$, we may assume WLOG that the
$J_\xi$ form a $\Delta$-system with some countable root $R$.
The $\PPP_R\,$-extension of $V$
still satisfies $\kappa = (2^{\aleph_0})^+$, so,
replacing $V$ by its $\PPP_R\,$-extension,
we may assume WLOG that $R = \emptyset$, so that the $J_\xi$ are 
disjoint.  We may also assume WLOG that the $|J_\xi|$
are all the same cardinal $\lambda \le \aleph_0$.
But now we may assume WLOG that $J_\xi = \{\xi\}$
and $\theta = \kappa$, since we may replace the $\PPP_i$
by the finite support product of $\lambda$ of the $\PPP_i$,
and simply discard the indices $i \in \theta \setminus \bigcup_\xi J_\xi$.
Now, each $\dot \CC_\xi$ is a $\PPP_\xi$-name.
WLOG, all the $\PPP_\xi$ are the same (since they are isomorphic),
so that whenever 
$\pi$ is a permutation of $\kappa$, it induces a 
natural automorphism $\widehat\pi$ of $\PPP$.
This automorphism also applies to the $\PPP$-names, so that
each $\widehat\pi (\dot \CC_\xi)$ is
a $\PPP_{\pi(\xi)}$-name.  But,
$|\PPP_\xi| \le 2^{\aleph_0}$, so
there are only $2^{\aleph_0}$ $\PPP_\xi$-names for countable
subsets of $\PP(\omega)$, so WLOG, we may assume that
$\widehat\pi (\dot \CC_\xi)$ is always the name $\dot \CC_{\pi(\xi)}$.

Still in $V$, we have a name
$\dot T$ for the tree $T( \langle \dot \CC_\xi : \xi < \kappa \rangle)$.
We also have a name $\dot\rho$ for the rank function.
Whenever $\vec \xi \in \kappa^{< \omega}$ is
a sequence of distinct ordinals, let
\[
S_{\vec \xi} =
\{\sigma : \exists p \in \PPP \, [ p \Vdash
[\exists \vec H \, [(\vec \xi, \vec H) \in \dot T
\ \&\ \dot\rho(\vec\xi,\vec H) = \sigma]]]\}\ \ .
\]
Then each $|S_{\vec\xi}| \le \aleph_0$
(since $\dot \CC_\xi$ is forced to be countable), and
$S_{\vec\xi}$ only depends on the length of $\vec\xi$
(using the automorphisms; note that in Definition \ref{def-tree},
the tree order does not depend on the ordering of the ordinal $\gamma$). 
Thus, if we set $S = \bigcup\{S_{\vec\xi} : \vec \xi \in \kappa^{< \omega} \}$,
then $|S| \le \aleph_0$.
Since $S$ must also be an initial segment of
the ordinals, it is forced by $\one$ that
$\rank(\dot T) < \omega_1$.

Now, in the random real case, we may assume (in $V$) that
the elements $p \in \PPP$ are the Baire subsets of $2^\theta$
of positive measure.  The order $\le$ is just $\subseteq$,
and $\one = 2^\theta$. 
If $J \subseteq \theta$, we 
let $\PPP_J$  be the set of elements $p \in \PPP$
such that $J$ is a support of $p$ (i.e., 
$\forall f,g \in 2^\theta\, [f\res J = g\res J \to
[f \in B \leftrightarrow g \in B]]$).
Every Baire subset of $2^\theta$ has a countable support, so
that each $\dot \CC_\xi$  is a $\PPP_{J}$-name for
some countable $J$.
The same countable $\Delta$-system argument lets us
assume that WLOG, each 
$\dot \CC_\xi$  is a $\PPP_{J_\xi}$-name, where
$J_\xi = \{\omega \cdot \xi + n : n \in \omega\}$.
Now, $\widehat\pi$ is the automorphism of $\PPP$ which
arises from permuting coordinate $\omega \cdot \xi + n$ to
coordinate $\omega \cdot \pi(\xi) + n$.
The proof that $\one \Vdash \rank(\dot T) < \omega_1$ is the same
as before.
\end{proof}

The following theorem yields $(B)$ in the extension when the ground model
satisfies CH:

\begin{theorem}
\label{thm-main}
In the ground model $V$, set $\kappa = (2^{\aleph_0})^+$.
Let $\PPP$ be ccc and be  either a random real extension
or of the form $\prod^{fin}_{i \in \theta} \PPP_i$,
where all the $\PPP_i$ are isomorphic and $|\PPP_i| \le 2^{\aleph_0}$.

Then in the generic extension $V[G]$, the following holds:
Suppose $X$, $\gamma$, and $E_\xi$, for $\xi < \gamma$,
satisfy $(1)$ and $(2)$ of Lemma \ref{lemma-bigrank},
and in addition all $|E_\xi| = \aleph_0$.  Then $\gamma < \kappa$.
In particular, whenever $X$ is compact scattered,
$|\{\alpha < \hgt(X) : |I(X^{(\alpha)})| = \aleph_0\}| <  \kappa$.
\end{theorem}
\begin{proof}
The fact that $\gamma < \kappa$ is immediate
from Lemmas \ref{lemma-omega} and \ref{lemma-smallrank}.
For the ``in particular'':  If that failed,
we could let $\langle \alpha_\xi : \xi < \kappa \rangle$
be an increasing sequence of ordinals with each
$|I(X^{(\alpha_\xi)})| = \aleph_0$; then setting  $\gamma = \kappa$ and
$E_\xi = I(X^{(\alpha_\xi)})$, we would have a contradiction.
\end{proof}

\begin{acknowledgement}
We would like to thank Istv\'an Juh\'asz for
some useful comments on the original draft of this paper.
\end{acknowledgement}

\end{document}